\documentclass[12pt]{amsart}
\usepackage{amsaddr}

\usepackage{amssymb}
\usepackage{mathabx}
\usepackage{bbm}
\usepackage{amsthm}
\usepackage{dsfont}
\usepackage{amsmath}
\usepackage{color,graphics,srcltx}
\usepackage{bm}
\usepackage{easybmat}
\usepackage{multirow,bigdelim}
\usepackage{enumerate}
\usepackage{graphicx}

\newtheorem{proposition}{Proposition}
\newtheorem{example}[proposition]{Example}
\newtheorem{lemma}[proposition]{Lemma}

\newtheorem{theorem}[proposition]{Theorem}

\newtheorem{remark}[proposition]{Remark}

\newcommand{\II}{\mathds{I}}

\begin{document}

\title{Marshall-Olkin copulas revisited${}^{1}$}
\thanks{${}^{1}$All authors acknowledge the support of ARIS (Slovenian Research and Innovation Agency) research core
funding No. P1-0448. The work of the second author was also supported by Croatian Science Foundation grant no. 2277.}

%\author[D.~Kokol Bukov\v{s}ek]{Damjana Kokol Bukov\v{s}ek} \address{University of Ljubljana, School of Economics and Business, and Institute of Mathematics, Physics and Mechanics, Ljubljana, Slovenia} \email{damjana.kokol.bukovsek@ef.uni-lj.si}

\author[T.~Ko\v{s}ir]{Toma\v{z} Ko\v{s}ir}
\address{University of Ljubljana, Faculty of Mathematics and Physics%, and Institute of Mathematics, Physics and Mechanics
, Ljubljana, Slovenia}
\email{tomaz.kosir@fmf.uni-lj.si}

\author[P.~Lazi\'{c}]{Petra Lazi\'{c}} \address{University of Ljubljana, Faculty of Mathematics and Physics, Ljubljana, Slovenia, and University of Zagreb, Faculty of Science, Department of Mathematics, Zagreb, Croatia}
\email{Petra.Lazic@fmf.uni-lj.si}

\author[M.~Omladi\v{c}]{Matja\v{z} Omladi\v{c}}
\address{University of Ljubljana, Faculty of Mathematics and Physics, and Institute of Mathematics, Physics and Mechanics,  Ljubljana, Slovenia}
\email{matjaz@omladic.net}

\begin{abstract}
    Almost seventy years old Marshall-Olkin copulas, then wider Marshall copulas, and finally even wider shock model (SM) copulas constitute a substantial part of nowadays copula theory due to numerous applications. Recently, Christian Genest with some coauthors introduced a new stochastic model for a special subclass of SM copulas which gives not only a new angle on these copulas but also widens the range of applications. In this paper we extend this type of stochastic model to all known subclasses of SM copulas. We also introduce a novel class of SM copulas and extend the new stochastic model to this class as well.
\end{abstract}

%\thanks{DKB\&TK\&BM\&MO acknowledge financial support from the Slovenian Research Agency (research core funding No. P1-0222).}

\maketitle

\section{ Introduction }\label{sec:intro}

Copula is simply a joint distribution with uniform margins, but when we insert arbitrary univariate distributions as margins  into it, we can get any bivariate distribution. This seminal 1959 result of Sklar \cite{Skla} has made them the most important tool of dependence modeling \cite{Joe,Nels,DuSe}. Most of what we do in this paper has an extension to its multivariate version, but we will restrict ourselves to the bivariate case to make intuitively clearer the stochastic model that we are introducing.

If one wants to build a stochastic model for describing the dependence among two (or more) lifetimes, i.e. positive random variables, one comes to shock models. In engineering applications, joint models of lifetimes may serve to estimate the expected lifetime of a system composed by several components. In a related situation like portfolio credit risk, instead, the lifetimes may have the interpretation of time-to-default of firms, or generally financial entities, while a stochastic model may estimate the price/risk of a related derivative contract (e.g. CDO). In both cases, it is of interest to estimate the probability of the occurrence of a joint default. This approach is a modern practitioners' view on what was started analytically in the 1967 key paper of Marshall and Olkin \cite{MaOl} and was extended later by Marshall in \cite{Mars}. (Perhaps we should be aware that in the original paper \cite{MaOl} the word copula was not even mentioned.)
An important portion of copula theory today (here is some recent publications \cite{DuGiMa,DuGiMa1,DuKoMeSe,GoKoMu,LiPe,Muli}) 
studies or is at least based on these two papers. The first one allows only exponential shocks, i.e. lifetimes, while the second one allows for arbitrary distributions. In 2016 Omladi\v c and Ru\v zi\' c \cite{OmRu} propose to study Marshall type of models which result in slightly different kind of copulas they call maxmin copulas. That paper started a series of papers that introduced some classes of copulas obtained in similar ways, including RMM (reflected maxmin) copulas, proposed by Ko\v sir and Omladi\v c in 2020 \cite{KoOm}, cf.\ also \cite{DuOmOrRu,KaTr,KoBuKoMoOm1,KoBuKoMoOm3,KoOm0}. 
Many authors are calling copulas of all these types \emph{Marshall-Olkin copulas}, but most of them keep this term only for the type introduced in \cite{MaOl}, i.e. the ones modelled by exponential shocks, and call the general class \emph{Marshall copulas}. We will call the whole class \emph{shock model copulas} and add the term \emph{exponential} respectively \emph{general} if our shocks will be distributed exponentially resp. according to more general distributions. %Since most of our copulas will use exponential distributions we will in this case omit the word exponential.

The main goal of this paper is twofold. First, we follow a new stochastic model that gives rise to the usual shock model copulas and has some advantages over the original one. This kind of model was first introduced by Genest et al. in 2018 \cite{GeMeSc} and later by Bentoumi et al. \cite{BeKtMe} and by El Ktaibi et al. \cite{KtBeSoMe}. They are developing the model only for very special classes of exponential models. We show that this approach extends to all known cases of shock model copulas. Moreover, we introduce a novel class we call the \emph{survival RMM copula (SMM for short)}. The proposal of this class is the second goal of this paper. Of course, we  extend the new stochastic model to this class as well.
%We provide a natural setup of shocks for which the dependence is modelled by reflected maxmin copulas, RMM for short (see \cite{KoOm,KoOm0,KaTr}). They are the reflected family of the maxmin copulas introduced by Omladi\v{c} and Ru\v{z}i\'{c} (see \cite{OmRu,DuOmOrRu}). In perhaps the most important special case when the shocks are all distributed exponentially we recover the model introduced by El Ktaibi et al.\ (see \cite{KtBeSoMe,BeKtMe}).

The paper is organized as follows: Section 2 provides preliminaries on the Marshall, Maxmin and RMM copulas. Section 3 exhibits the new stochastic model for the Marshall case. Section 4 presents the new model for the RMM case and Section 5 introduces the new class of SMM copulas. Finally, Section 6 gives the new stochastic model for the SMM copulas.

\section{Preliminaries on Marshall, maxmin, and RMM copulas}

A function $C: [0,1]\times[0,1] \to [0,1]$ is called a \emph{copula} if conditions \emph{(a)} and \emph{(b)} below are fulfilled. %A \emph{copula} may be seen as a quasi-copula $C$ such that $C$-volume of every rectangle is positive.
\emph{\begin{enumerate}[(a)]
  \item $C$ is grounded and 1 is the {neutral element} of $C$;
  \item $C$-volume of every rectangle is positive.
\end{enumerate}}
\noindent
Here, condition \emph{(a)} means that $C(x, 0) = C(0,x)=0$ for all $x \in [0, 1]$ and $C(x, 1) = C(1,x)= x$ for all $x \in [0, 1]$. Condition \emph{(b)} means that
\emph{\begin{enumerate}[(P)]
        \item $V_C([x_1, x_2] \times [y_1,y_2]) = C(x_2,y_2) - C(x_1, y_2) - C(x_2,y_1)+ C(x_1, y_1) \geq 0$  for all $x_1,x_2,y_1,y_2 \in [0, 1]$  such that $x_1\leq x_2$ and $y_1\leq y_2$.
      \end{enumerate}}
\noindent A set of the form $[x_1, x_2] \times [y_1,y_2]\subseteq[0,1]\times[0,1]$, $x_1\leq x_2,y_1\leq y_2$, is called a \emph{rectangle} and condition \emph{(P)} may be seen as the definition of its $C$-volume denoted by $V_C$, together with the requirement that it is positive.\\[1mm]

The copulas presented here are sometimes called bivariate copulas as opposed to more general multivariate copulas, where the bivariate domain $(x,y)\in[0,1]^2$ is replaced by $(x_1,x_2,\ldots,x_n)\in[0,1]^n$ for any integer $n\geq2$. However, we will not study copulas i\emph{}n this generality. It turns out that for any copula $C$ we have 
\[
    W(x,y)\leqslant C(x,y)\leqslant M(x,y),
\]
where 
\[
    W(x,y)=\max\{0,x+y-1\}\quad\mbox{respectively}\quad M(x,y)=\min\{x,y\}
\]
are copulas called respectively Fr\' echet-Hoeffding lower bound and upper bound of all bivariate copulas.

For a bivariate distribution $F(x,y)$ of a random vector $(X,Y)$, it is easy to find the distributions of its components. For simplicity we assume that the domains of $X$ and $Y$ are contained in the interval $[0,1]$. Then
\[
    F_X(x)=F(x,1)\quad\mbox{and}\quad F_Y(y)=F(1,y)
\]
are the respective distributions of the components. We call them the \emph{marginal distributions}. (For instance, copulas may be seen as bivariate distributions whose marginal distributions are both uniform on the interval $[0,1]$.) In his seminal 1959 paper \cite{Skla} A. Sklar defined copulas and proved the following theorem. (Actually, he proved it in the multivariate setting.)

\noindent{\bf Sklar Theorem.}\emph{ For any bivariate distribution $F(x,y)$ there exists a copula $C(x,y)$ such that
\[
    F(x,y)=C(F_X(x),F_Y(y))\quad\mbox{where}\quad F_X,F_Y
\]
are the corresponding marginal distributions. Conversely, if $C(x,y)$ is any copula and $F_X,F_Y$ are any univariate distributions, then
\[
    H(x,y)=C(F_X(x),F_Y(y))
\]
is a bivariate distribution whose copula is $C(x,y)$.}

We refer to \cite{DuSe,Nels} for further details on general theory on copulas.\\[1mm]

Marshall \cite{Mars} showed that shock models governed by
\[
    U=\max\{X,Z\},\quad V=\max\{Y,Z\},
\]
where $X$, $Y$, and $Z$ are independent random variables, are described by copulas of the form 
\[
    C(u,v)= \min\{ug(v),f(u)v\},
\]
where generating functions $f$ and $g$ satisfy the following conditions (here and in what follows we are using notation $\II$ for the closed interval $[0,1]$)
\begin{enumerate}[(M1)]
  \item $f(0)=0,\ f(1)=1$,
  \item $f$ is nondecreasing on $\II$,
  \item function $f^*(u)=\dfrac{f(u)}{u}$ is nonincreasing on $(0,1]$.
\end{enumerate}
Observe that (M3) implies that $f(u)\geqslant u$ for all $u\in\II$. 
The joint distribution function $H$ of $(U,V)$ is then given by
\[
    H(x,y)= C(F_U(x),F_V(y)),
\]
where $F_U$ and $F_V$ are the marginal distribution functions of $(U,V)$, i.e., distributions of $U$ and $V$ respectively given by 
\[
    F_U(x)=F_X(x)F_Z(x)\quad\mbox{and}\quad F_V(x)=F_Y(x)F_Z(x).
\]

Omladi\v{c} and Ru\v{z}i\'{c} \cite{OmRu} introduced a new class of copulas modeling shocks given by 
\[
    U=\max\{X,Z\},\quad V=\min\{Y,Z\}.
\]
where $X$, $Y$, and $Z$ are independent random variables. A copula for $(U,V)$ is then of the form 
\begin{equation}\label{eq:omru}
  C(u,v)= \min\{u,\phi(u)(v-\psi(v))+u\psi(v)\},
\end{equation}
where generating function $\phi$ satisfies conditions (M1)-(M3) with $\phi=f$, while $\psi$ has the following properties
\begin{enumerate}[(F1)]
  \item $\psi(0)=0,\ \psi(1)=1$,
  \item $\psi$ is nondecreasing on $\II$,
  \item function $\psi_*(v)=\dfrac{1-\psi(v)}{v-\psi(v)}$ is nonincreasing on $[0,1)$. Here, $\psi_*(v)=\infty$ if $\psi(v)=v$.
\end{enumerate}
Functions of the form \eqref{eq:omru} are called \emph{maxmin copulas}.
 
Ko\v{s}ir and Omladi\v{c} \cite{KoOm} studied reflected version of the maxmin copulas that are called \emph{reflected maxmin copulas} (\emph{RMM} for short). Their generators $f$ and $g$ satisfy the following conditions. (Here, we introduce additionally functions $f^*(u)=\dfrac{f(u)}{u}$, $g^*(u)=\dfrac{g(u)}{u}$,  $\widehat{f}(u)=f(u)+u$, and $\widehat{g}(u)=g(u)+u$  in order to make these conditions easier to state and understand.)
\begin{enumerate}[(G1)]
  \item $f(0)=g(0)=0,\ f(1)=g(1)=0$, and $f^*(0)=g^*(0)=0$;
  \item functions $\widehat{f}$ and $\widehat{g}$ are nondecreasing on $\II$;
  \item functions $f^*$ and $g^*$ are nonincreasing on $\II$.
\end{enumerate}
An RMM copula is then of the form
\[
    C(u,v)=\max\{0, uv-f(u)g(v)\},
\]

\section{Shock model for Marshall copulas}\label{sec:setup}

We first consider a slight variation of Marshall model \cite[Proposition 3.2]{Mars} that still leads to the same family of copulas. We consider four random variables $X$, $Y$, $Z_1$, and $Z_2$ representing shocks. The two triples $X,Y,Z_1$ and  $X,Y,Z_2$ are independent, while $Z_1$ and $Z_2$ form a comonotonic pair, so that their dependence is modeled by the Fr\' echet-Hoeffding upper bound $M(u,v)=\min\{ u,v\}$. Denote by $F_X$ and $F_Y$ the distribution functions of $X$ and $Y$ respectively, and by  $G_i$ the distribution functions of $Z_i$ for $i=1,2$. We construct two random variables
\[
    U=\max\{X,Z_1\},\ \ V=\max\{Y,Z_2\},
\]
whose distributions are clearly $F_U(x)=F_X(x)G_1(x)$ and  $F_V(y)=F_Y(y)G_2(y)$. Before giving the copula connecting $U$ and $V$ we introduce some further notation.

Let $F$ be a distribution function. For $u\in\II$ we have %denote by $F^{-1}$ its quasi-inverse defined by
\[
    F^{-1}(u)=\inf\{x\in\mathds{R}\,:\,F(x)\geqslant u\}.
\]
We assume that infimum of an empty set is equal to $+\infty$ and infimum of $\mathds{R}$ is equal to $-\infty$. Note that $F^{-1}$ is nondecreasing. We have $F(F^{-1}(u))\geqslant u$ for all $u\in\II$ and  $F^{-1}(F(x))\leqslant x$ for all $x\in\mathds{R}$. Following the notation of \cite{OmRu} we denote further by $f(x-)$ the left limit of function $f$ if it exists. (Note that all our functions will be monotone, so that our left limits will always exist at $x$.) For $u\not\in\mathrm{im}F\cup\{0,1\}$ we write $\overline{u}=F(F^{-1}(u))$ and $\underline{u}=F(F^{-1}(u)-)$. Then $\overline{u}\in\mathrm{im}F$ and either $\underline{u}\in\mathrm{im}F$ or $\underline{u}\not\in\mathrm{im}F$ and $(\underline{u}-\varepsilon,\underline{u}) \cap\mathrm{im}F\neq\emptyset$ for every small enough $\varepsilon>0$.

%Now, suppose that $X$ and $Y$ are independent random variables, while $Z_1$ and $Z_2$ are counter-monotonic random variables, that is   \[    Z_1=G_1^{-1}(Z)\quad \mbox{and}\quad Z_2=G_2^{-1}(I-Z),\]where $Z$ is a random variable distributed uniformly on $\II$, and such that $X,Y,Z$ are independent. We want to model the dependence between \[     U=\max\{X,Z_1\},\ \ V=\max\{Y,Z_2\}.\]We denote by $F_U$ and $F_V$ the distribution functions of $U$ and $V$ respectively.

Let us introduce generating functions for our model by
\[
    \varphi(u)=\begin{cases}
                 0, & \mbox{if } u=0 \\
                 1, & \mbox{if } u=1 \\
                 F_X(F_U^{-1}(u)), & \mbox{if } u\in\mathrm{im}F_U\setminus\{0,1\} \\
                 \dfrac{\varphi(\overline{u})-\varphi({\underline{u}}-)}{\overline{u}-\underline{u}} (u-\underline{u})+\varphi(\underline{u}-) & \mbox{otherwise};
               \end{cases}
\]
and
\[
    \psi(v)=\begin{cases}
                 0, & \mbox{if } v=0 \\
                 1, & \mbox{if } v=1 \\
                 F_Y(F_V^{-1}(v)), & \mbox{if } v\in\mathrm{im}F_V\setminus\{0,1\} \\
                 \dfrac{\psi(\overline{v})-\psi({\underline{v}}-)}{\overline{v}-\underline{v}} ({v}-{\underline{v}})+\psi({\underline{v}}-) & \mbox{otherwise}.
               \end{cases}
\]
A careful reader will have noticed that the first three lines of each definition are the obvious natural choice while the fourth one is a linear interpolation on the intervals left undefined. %The generating functions are then \[     f(u)=\widehat{f}(u)-u\quad\mbox{and}\quad g(v)=\widehat{g}(v)-v.\]
Let us give the proof of the following proposition for the sake of completeness.

\begin{proposition}\label{pr:mars}
  The copula of the pair $(U,V)$ is a Marshall copula
  \[
    C(u,v)=\min\{u\psi(v),v\varphi(u)\}.
  \]
\end{proposition}

\begin{proof}
  Denote by $H$ the joint distribution function of the random vector $(U,V)$ and compute
  \[
    \begin{split}
       H(x,y) & =\mathds{P}[U\leqslant x,V\leqslant y]=\mathds{P}[X\leqslant x,Z_1\leqslant x,Y\leqslant y,Z_2\leqslant y] \\
         & =F_X(x)F_Y(y)\mathds{P}[Z_1\leqslant x,Z_2\leqslant y]=F_X(x)F_Y(y)\min\{G_1(x),G_2(y)\}\\
         & =\min\{F_U(x)\psi(F_V(y)),F_V(y)\varphi(F_U(x))\}\\
         & =C(F_U(x),F_V(y)),
    \end{split}
  \]
  where $C$ is the Marshall copula with generating functions $\varphi$ and $\psi$.
\end{proof}

We now state and prove two results that exhibit a modified version of Marshall's Theorem \cite[Proposition~3.2]{Mars}. Our results are based on Lemma 5.6 and Theorem 5.9 of \cite{Ruzic}.

\begin{lemma}\label{le:mars}
  Suppose that $U$ and $V$ are two random variables with distribution functions $F_U$ and $F_V$, and with the joint distribution function $H$. Then the following two statements sre equivalent. 
  \begin{enumerate}[(a)]
    \item There are two triples of independent random variables $X,Y,Z_1$ and $X,Y,Z_2$, where $Z_1$ and $Z_2$ are comonotonic, such that
        \[
            U=\max\{X,Z_1\}\quad\mbox{and}\quad V=\max\{Y,Z_2\}.
        \]
    \item There are distribution functions $F_X, F_Y, G_1,$ and $G_2$ such that
        \[
            \begin{split}
               H(x,y) & =F_X(x)F_Y(y)\min\{G_1(x),G_2(y)\} \\
                 & =F_X(x)F_Y(y)M(G_1(x),G_2(y))
            \end{split}
        \]
  \end{enumerate}
  If any of (a) or (b) holds, then we also have
  \[
    F_U(x)=F_X(x)G_1(x),\quad F_V(y)=F_Y(y)G_2(y),
  \]\[
    F_U\leqslant\min\{F_X,G_1\},\quad\mbox{and}\quad F_V\leqslant\min\{F_Y,G_2\}.
  \]
\end{lemma}

\begin{proof}
  Suppose \emph{(a)} holds. Then, we see as in the proof of Proposition \ref{pr:mars} that
        \[
            H(x,y)=F_X(x)F_Y(y)\min\{G_1(x),G_2(y)\}.
        \]
  Now, assume \emph{(b)}. Then, $H(x,y)$ is the joint distribution function of random variables $U=\max\{X,Z_1\}$ and $V=\max\{Y,Z_2\}$ as above. Since the random vector is uniquely determined by its joint distribution function, \emph{(a)} follows.
\end{proof}

\begin{remark}\label{rem}
  Lemma \ref{le:mars} remains correct if we replace the word ``comonotonic'' in (a) with the word ``countermonotonic'' and the formula in (b) with
        \[
               H(x,y) =F_X(x)F_Y(y)W(G_1(x),G_2(y))  
        \]
  The proof follows in a similar way as above. \\[1mm]
\end{remark}

\begin{theorem}\label{thm:3}
  Let $C_{\varphi,\psi}$ be a Marshall copula, and let $F_U$ and $F_V$ be two distribution functions that satisfy the following assumptions:
  \begin{enumerate}[(a)]
    \item There is an increasing %c\`adl\`ag 
    function $\chi:\mathds{R}\rightarrow\mathds{R}$ such that if $x\in\mathds{R}$ is such that $F_U(\chi(x))>0$ and $F_V(x)>0$, then
        \[
            \varphi^*(F_U(\chi(x)))=\psi^*(F_V(x)).
        \]
    \item Function $\varphi$ is either continuous at $0$ or \[x_U=\inf\{x\in\mathds{R},F_U(x)>0\}>-\infty\] and $F_U$ is not continuous at $x_U$. Function $\psi$ is either continuous at $0$ or \[x_V=\inf\{x\in\mathds{R}, F_V(x)>0\}>-\infty\] and $F_V$ is not continuous at $x_V$.
    \item Function $\varphi$ is such that  either $\varphi^*(0+)=\infty$ or there is $x\in\mathds{R}$ such that $F_U(x)=0$. Function $\psi$ is such that either $\psi^*(0+)=\infty$ or there is $x\in\mathds{R}$ such that $F_V(x)=0$. 
  \end{enumerate}
  Then there are two triples of independent random variables $X,Y,Z_1$, and $X,Y,Z_2$, where $Z_1$ and $Z_2$ are comonotonic, such that $H(x,y)= C_{\varphi.\psi}(F_U(x),F_V(y))$ is the joint distribution function of 
  \[
    U=\max\{X,Z_1\}\quad\mbox{and}\quad V=\max\{Y,Z_2\}.
  \]
\end{theorem}

\begin{proof}
  We define $F_X(x)=\varphi(F_U(x))$ and $F_Y(y)=\psi(F_V(y))$ for all $x,y\in\mathds{R}$. Assumptions \emph{(b)} and \emph{(c)} imply that $F_X$ and $F_Y$ are distribution functions. Next we define
  \[
    G_2(x)=\begin{cases}
             0, & \mbox{if } F_U(\chi(x))=F_V(x)=0, \\
             \dfrac{F_V(x)}{\psi(F_V(x))}, & \mbox{if } F_U(\chi(x))=0\ \mbox{and}\ F_V(x)>0, \\
             \dfrac{F_U(\chi(x))}{\varphi(F_U(\chi(x)))}, & \mbox{otherwise}.
           \end{cases}
  \]
  Note that
  \[
    \dfrac{F_V(x)}{\psi(F_V(x))}=\dfrac{F_U(\chi(x))}{\varphi(F_U(\chi(x)))}
  \]
  whenever both  $F_U(\chi(x))>0\,\ {and}\ F_V(x)>0$ by assumption \emph{(a)}. We also take $G_1(x)=G_2(\chi^{-1}(x))$ for all $x\in\mathds{R}$. Here, $\chi^{-1}$ is the inverse of $\chi$ which is an increasing function since $\chi$ is increasing. Functions $G_1$ and $G_2$ are well defined because $\varphi(x)\geqslant0$ and $\psi(y)\geqslant0$ for all $x,y\in\II$.
  Assumptions \emph{(a)--(c)} imply that $G_1$ and $G_2$ are distribution functions. 
  % Furthermore, we have
  If $F_U(x)=0$ then
  \[
    F_X(x)G_1(x)=\varphi(F_U(x))G_1(x)=0=F_U(x).
  \]
  If $F_U(x)>0$, then
  \[
    F_X(x)G_1(x)=\varphi(F_U(x))G_2(\chi^{-1}(x))=\dfrac{\varphi(F_U(x))F_U(x)}{\varphi(F_U(x))}=F_U(x).
  \]
  Therefore we have $F_U(x)=F_X(x)G_1(x)$  for all $x\in\mathds{R}$. 
  If $F_V(y)=0$ then
  \[
    F_Y(y)G_2(y)=0=F_V(y).
  \]
  Otherwise, if $F_V(y)>0$ then
  \[
    F_Y(y)G_2(y)=\dfrac{\psi(F_V(y))F_V(y)}{\psi(F_V(y))} =F_V(y). %-*(F_V(y))G_1(\chi(y))=F_V(y)\varphi^*(F_U(\chi(y)))G_1(\chi(y))=F_V(y).
  \]
  Hence we have $F_V(y)=F_Y(y)G_2(y$ for all $y\in\mathds{R}$.
  
  Finally, the definitions of $H, C_{\varphi,\psi}, F_X, F_Y, G_1,$ and $G_2$ imply
  \[
    \begin{split}
       H(x,y) & = C_{\varphi,\psi}(F_U(x),F_V(y))=\min\{F_U(x)\psi(F_V(y)),\varphi(F_U(x))F_V(y)\} \\
         & =\min\{F_U(x)F_Y(y),F_V(y)F_X(x)\}= F_X(x)F_Y(y)\min\{G_1(x),G_2(y)\}.\\
    \end{split}
  \]
  The implication \emph{(a)}$\Rightarrow$\emph{(b)} of Lemma \ref{le:mars} then completes the proof.
\end{proof}

\section{Shock model for RMM copulas}\label{sec:setup2}

Next we present a shock model setup first proposed in Genest et al.\ \cite{GeMeSc} (cf.\ also \cite{BeKtMe,KtBeSoMe}). 
We consider two triples $(X,Y,Z_1)$ and $(X,Y,Z_2)$ of independent  random variables. This time we assume that $Z_1$ and $Z_2$ are a countermonotonic pair, so that their dependence is measured by the Fr\' echet-Hoeffding lower bound $W(u,v) =\max\{0,u+v-1\}$. Again we construct two random variables
\[
    U=\max\{X,Z_1\},\ \ V=\max\{Y,Z_2\},
\]
whose distributions are clearly $F_U(x)=F_X(x)G_1(x)$ and  $F_V(y)=F_Y(y)G_2(y)$ as above. 
Let us introduce auxiliary generating functions for our model by
\[
    \widehat{f}(u)=\begin{cases}
                 0, & \mbox{if } u=0 \\
                 1, & \mbox{if } u=1 \\
                 F_X(F_U^{-1}(u)), & \mbox{if } u\in\mathrm{im}F_U\setminus\{0,1\} \\
                 \dfrac{\widehat{f}(\overline{u})-\widehat{f}({\underline{u}}-)}{\overline{u}-\underline{u}} (u-\underline{u})+\widehat{f}(\underline{u}-) & \mbox{otherwise};
               \end{cases}
\]
and
\[
    \widehat{g}(v)=\begin{cases}
                 0, & \mbox{if } v=0 \\
                 1, & \mbox{if } v=1 \\
                 F_Y(F_V^{-1}(v)), & \mbox{if } v\in\mathrm{im}F_V\setminus\{0,1\} \\
                 \dfrac{\widehat{g}(\overline{v})-\widehat{g}({\underline{v}}-)}{\overline{v}-\underline{v}} ({v}-{\underline{v}})+\widehat{g}({\underline{v}}-) & \mbox{otherwise}.
               \end{cases}
\]
Note that $F_U^{-1}(0) = - \infty$, so $F_X (F_U^{-1}(0)) = 0$ independently of the point 0 being in im$\,F_U$ or not, and similarly for $\widehat{g}$. In the same way, if 1 is not in im$\,F_U$, we get that $F_X (F_U^{-1}(1)) = F_X(+\infty) = 1$, and if 1 is in im$\,F_U$, then also $F_X (F_U^{-1}(1)) = 1$.
We write also $f(u)=\widehat{f}(u)-u$ and $g(v)=\widehat{g}(v)-v$.
%Before giving the copula connecting $U$ and $V$ we introduce some further notation.

\begin{theorem}\label{thm:5}
  The copula of the pair $(U,V)$ is an RMM copula with generating functions $f$ and $g$
  \[
    C(u,v)=\max\{0,uv- f(u)g(v)\}.
  \]
\end{theorem}

\begin{proof}
  Since $X$ and $Z_1$, respectively $Y$ and $Z_2$, are independent, we have
  \begin{equation}\label{eq:rmm}
    F_U(u)=F_X(u)G_1(u),\quad\mbox{respectively}\quad F_V(u)=F_Y(u)G_2(u).
  \end{equation}
  The joint distribution function $H$ of the random vector $(U,V)$ is then
  \[
    \begin{split}
        H(u,v) & =\mathds{P}[U\leqslant u, V\leqslant v] \\
         & =\mathds{P}[X\leqslant u,Y\leqslant v,G_1^{-1}(Z)\leqslant u,G_2^{-1}(1-Z)\leqslant v ].
    \end{split}
  \]
  Since $X,Y$, and $Z$ are independent, we obtain
  \[
    \begin{split}
       H(u,v) & =F_X(u)F_Y(v)\mathds{P}[1-G_2(v)\leqslant Z\leqslant G_1(u)] \\
         & =F_X(u)F_Y(v) \max \{G_1(u)-(1-G_2(v)),0\} \\
         & =F_X(u)F_Y(v)W(G_1(u),G_2(v)),
    \end{split}    
  \]
  where in the second line we have used the fact that $Z$ is distributed uniformly on $\II$ and then the definition of the Fr\'{e}chet-Hoeffding lower bound $W$. % From now on the notation $W$ for the Fr\'{e}chet-Hoeffdung lower bound will be used throughout.
  The copula of $(U,V)$ is then given by 
  \[
  \begin{split}
     C(u,v) & =H(F_U^{-1}(u),F_V^{-1}(v)) \\
       & =F_X(F_U^{-1}(u))F_Y(F_V^{-1}(v))W(G_1(F_U^{-1}(u)),G_2(F_V^{-1}(v))).
  \end{split}    
  \]
  We use \eqref{eq:rmm} and the definition of the generating functions $f$ and $g$ to conclude
  \[
    C(u,v)=\max\{0,uv- f(u)f(v)\}.
  \]
  It remains to show that functions $f$ and $g$ satisfy conditions (G1)-(G3) for generators of RMM copulas. 
  
  (G1): By definition, $\widehat{f}(0) = 0$ and $\widehat{f} (1) = 1$.
  %First we see that $\widehat{f}(0)=F_X(F_U^{-1}(0))=F_X(-\infty)=0$. Then, denote $u_0=\inf\{u\in\mathds{R}\,:\,F_U(u)=1\}$ to get $F_U(u)=1$ for all $u\geqslant u_0$ and so $F_X(u)=G_1(u)=1$ as well. Also, $\widehat{f}(1)=F_X(F_U^{-1}(u_0))=F_X(1)=1$.
  
  (G2): Consider function $\widehat{f}(u)=f(u)+u=F_X(F_U^{-1}(u))$; since distribution functions $F_X$ and $F_U$ are nondecreasing, so that $F_U^{-1}$ is also nondecreasing, and since composite of two nondecreasing functions is nondecreasing, the claim follows for $\widehat{f}$. The conclusion for $\widehat{g}$ is shown in a similar way.
  
  (G3): We want to show that function $f^*(u)=\dfrac{f(u)}{u}$ is nonincreasing on $(0,1]$. It suffices to show that $f^*(u)+1=\dfrac{\widehat{f}(u)}{u}$ is nonincreasing on $(0,1]$. Choose $u\in\mathrm{im}F_U$ to get (using \eqref{eq:rmm})
  \[
        \dfrac{\widehat{f}(u)}{u}=\dfrac{F_X(F_U^{-1}(u))}{F_U(F_U^{-1}(u))} =\dfrac{F_X(F_U^{-1}(u))}{F_X(F_U^{-1}(u))G_1(F_U^{-1}(u))} =\dfrac{1}{G_1(F_U^{-1}(u))}.
  \]
  This function has a nondecreasing function in its denominator, so that it is nonincreasing on $\mathrm{im} F_U$. Since we defined $\widehat{f}$ outside this region by linear interpolation, the claim follows. Using similar arguments we prove that $g^*$ is also nonincreasing.
\end{proof}

\begin{theorem}\label{thm8}
  Let $U$ and $V$ be two random variables with respective distribution functions $F_U$ and $F_V$ such that there exists an $x_0\in\mathds{R}$ with $F_U(x_0),F_V(x_0)\in(0,1)$. Furthermore, let the copula of the pair $(F_U,F_V)$ be an RMM copula with generating functions $f$ and $g$, i.e.,
  \[
    C(u,v)=\max\{0,uv-f(u)g(v)\}.
  \]
  Then there exist random variables $X,Y,Z_1,Z_2$ such that $X,Y,Z_i$ are independent for $i=1,2$, $Z_2$ is countermonotonic to $Z_1$, and function
  \[
    H(x,y)=C(F_U(x),F_V(y))
  \]
  is the joint distribution function of the variables 
  \[\max\{X,Z_1\}\quad\mbox{and}\quad \max\{Y,Z_2\}.\]
\end{theorem}

\begin{proof}
  For all $x,y\in\mathds{R}$ define $F_X(x)=\hat{f}(F_U(x))$ and $F_Y(y)=\hat{g}(F_V(y))$, where $\hat{f}$ and $\hat{g}$ are defined just before Theorem \ref{thm:5}. Moreover, let 
  \[
    G_1(x)=\begin{cases}
             \dfrac{g^*(1-F_V(x))}{1+g^*(1-F_V(x))}, & \mbox{if } F_U(x)=0 \\
             \dfrac{F_U(x)}{\hat{f}(F_U(x))}, & \mbox{otherwise},
           \end{cases}
  \]
  and
  \[
    G_2(y)=\begin{cases}
             \dfrac{f^*(1-F_U(y))}{1+f^*(1-F_U(y))}, & \mbox{if } F_V(y)=0 \\
             \dfrac{F_V(y)}{\hat{g}(F_V(y))}, & \mbox{otherwise}.
           \end{cases}
  \]
  Since $f$ and $g$ are nondecreasing, the same holds for $\hat{f}$ and $\hat{g}$, so that $F_X$ and $F_Y$ are nondecreasing and c\` adl\` ag. Continuity of $\hat{f}$ and $\hat{g}$ implies that
  \[
    \lim_{x\to-\infty} F_X(x)=0\quad\mbox{and}\quad\lim_{x\to+\infty} F_X(x)=1,
  \]
  and similarly for $F_Y$. We conclude that $F_X$ and $F_Y$ are distribution functions of some random variables we denote by $X$ and $Y$. The condition given in the theorem tells us that if $F(x)=0$ for some $x$, then $x<x_0$, so that $F_V(x)\leqslant F_V(x_0)<1$. Similarly, $F_V(x)=1$ implies $x_0<x$, so that $F_U(x)\geqslant F_U(x_0)>0$.
  
  Now, choose $x\in\mathds{R}$ with $F_U(x)>0$ and $F_V(x)<1$ to get 
  \[
     \dfrac{1}{G_1(x)}=  \dfrac{1+g^*(1-F_V(x))}{g^*(1-F_V(x))}= \dfrac{\hat{f}(F_U(x))}{F_U(x)}%\quad\mbox{and}\quad \dfrac{1}{G_2(x)}=  \dfrac{\hat{g}(F_V(y))}{F_V(y)}.
  \]
  and similarly for $G_2$. When $x$ goes to $-\infty$, we use the fact that $g^*$ is continuous at $0$ to get that $G_1$ approaches to $0$. When $x$ goes to $+\infty$ we use the fact that $\hat{f}$ is continuous at $1$ to get that $G_1$ approaches to $1$. Next we prove that $G_1$ is nondecreasing. If $F_U(x)>0$ then
  \[
     \dfrac{1}{G_1(x)}=   \dfrac{{f}(F_U(x))}{F_U(x)} +1
  \]
  which is nonincreasing by (G3). If $F_V(x)<1$ then
  \[
     \dfrac{1}{G_1(x)}=  \dfrac{1}{g^*(1-F_V(x))}+1
  \]
  which is nonincreasing by (G3). So, it remains to consider the case that $F_U(x)=0$ and $F_V(y)=1$ for some $x$ and $y$. Clearly, we have $x<y$ and the point $x_0$ given in the theorem belongs to the interval $(x,y)$ so that
  \[
    \begin{split}
       G_1(x) & =\dfrac{g^*(1-F_V(x))}{1+g^*(1-F_V(x))}\leqslant\dfrac{g^*(1-F_V(x_0))}{1+g^*(1-F_V(x_0))} \\
         &  = \dfrac{\hat{f}(F_U(x_0))}{F_U(x_0)}\leqslant\dfrac{\hat{f}(F_U(y))}{F_U(y)}=G_1(y).
    \end{split}
  \]
  The first one of the inequalities above holds because $F_V(x),F_V(x_0)<1$ and the second one because $F_U(x_0),F_U(y)>0$. Furthermore, since $\hat{f}$ and $g^*$ are continuous and there exists some $x_0$ with $F_V(x)<1$ and $F_U(x_0)>0$, it follows that $G_1$ is c\` adl\` ag. Similar considerations and conclusions apply to $G_2$. Finally, we have shown existence of random variables $Z_1$ and $Z_2$ whose distribution functions are respectively $G_1$ and $G_2$.
  
  We now want to study the connection of the joint distribution function $H$ of $U$ and $V$ to the copula
  \[
    C(u,v)=\max\{0,uv-f(u)g(v)\}=uv\max\{1-f^*(u)g^*(v)\}.
  \]
  Indeed, after introducing marginal distributions $F_U$ and $F_V$ into this formula, we get
  \[
    \begin{split}
        &C(F_U(x),F_V(y))= \\
        & = F_U(x)F_V(y) \max\left\{0,1-\left(\dfrac{\hat{f}(F_U(x))-F_U(x)}{F_U(x)}\right)
    \left(\dfrac{\hat{g}(F_V(y))-F_V(y)}{F_V(y)}\right)\right\} \\
         & = F_U(x)F_V(y) \max\left\{0,\dfrac{\hat{f}(F_U(x))}{F_U(x)}+\dfrac{\hat{g}(F_V(y))}{F_V(y)}-
         \dfrac{\hat{f}(F_U(x))}{F_U(x)}\dfrac{\hat{g}(F_V(y))}{F_V(y)}\right\}\\
         & = \hat{f}(F_U(x))\hat{g}(F_V(y))\max\left\{0,\dfrac{F_U(x)}{\hat{f}(F_U(x))}+ \dfrac{F_V(y)}{\hat{g}(F_V(y))} -1 \right\}\\
         & = F_X(x)F_Y(y)W(G_1(x),G_2(y)).
    \end{split}
    \]
    Choose independent random varables with respective distributions $F_X,F_Y$, and $G_1$. Define $Z_2= G_2^{-1}(1-G_1(Z_1))$ to get a countermonotonic random variable whose distribution equals $G_2$.
    
  The proof that $H(x,y)$ is the joint distribution function of random variables $\max\{X,Z_1\}$ and $\max\{Y,Z_2\}$ follows by Remark \ref{rem}.
  \end{proof}

\begin{example}\label{Ex:8}
  The so-called \emph{Eyraud-Farlie-Gumbel-Morgenstern (EFGM) distributions} have been considered by many authors (cf.\ \cite{SaPlKoSeMeKl}). We will follow the approach of Durante and Sempi \cite[Section 6.3]{DuSe}. 
  We want to present shocks $X,Y,Z_1$ and $Z_2$ such that the corresponding RMM copula is of the form
  \[
    C_a(u,v)=uv-a^2uv(1-u)(1-v),\quad\mbox{where}\quad 0\leqslant a\leqslant 1.
  \]
  So, $C_a$ is of the EFGM type. (We assume actually $0<a\leqslant 1$.)
\end{example}

\begin{proof}
  Observe that $\widehat{f}_a(t)=(a+1)t-at^2$ and that in the model $\widehat{f}_a(t)=F_X(F_U^{-1}(t))$. Choose
  \[
    F_X(x)=\begin{cases}
             0, & \mbox{if } x\leqslant0 \\
             x, & \mbox{if } 0\leqslant x\leqslant1 \\
             1, & \mbox{otherwise},
           \end{cases}
  \]
  i.e., $X\sim U(0,1)$. To obtain $F_U$ and then also $G_1$, we have to solve 
  \[
    (a+1)t-at^2=x
  \]
  for $t$. We obtain
  \[
    F_U(x)=\begin{cases}
             0, & \mbox{if } x\leqslant0 \\
             \dfrac{1}{2a}\left(a+1-\sqrt{(a+1)^2-4ax}\right), & \mbox{if } 0\leqslant x\leqslant1 \\
             1, & \mbox{otherwise}.
           \end{cases}
  \]
  It is easy to verify that $F_U$ is a distribution function of a continuous random variable whose density function is equal to
  \[
    p_U(x)=\begin{cases}
             \dfrac{1}{\sqrt{(a+1)^2-4ax}}, & \mbox{if } 0<x<1 \\
             0, & \mbox{otherwise}.
           \end{cases}
  \] 
  To obtain the distribution function of $G_1$ we use the relation $F_U(x))=F_X(x)G_1(x)$. Therefore, we get after a simplification that
  \[
    G_1(x)=\begin{cases}
             \dfrac{2}{(a+1)+\sqrt{(a+1)^2-4ax}}, & \mbox{if } x\leqslant1 \\
             1, & \mbox{otherwise}.
           \end{cases}
  \]
  With not much more effort we can also compute the density of this continuous random variable
  \[
    p_1(x)=\begin{cases}
             \dfrac{2a}{\left((a+1)^2-2ax\right)\sqrt{(a+1)^2-4ax}+(a+1)^3-4(a+1)ax}, & \mbox{if } x<1 \\
             0, & \mbox{otherwise}.
           \end{cases}
  \]
  Since $C_a(u,v)$ is symmetric we take $F_X(x)=F_Y(x)$ and $G_1(x)=G_2(x)$, while $Z_1$ and $Z_2$ are related via copula $W(u,v)$. 
  \end{proof}

\begin{center}
\begin{figure}[h!]
      \caption{%Images corresponding to\\ 
      \textsl{Up -- Example \ref{Ex:8}:} Copula for the value of the parameter $a = 0.95$.\\ \textsl{Down -- Example \ref{Ex:9}:} Copulas for $\alpha = \beta = 0.1$ (left) and $\alpha = 0.4, \beta = 0.9$ (right). }\label{fig-bimod}
\includegraphics[width=5cm]{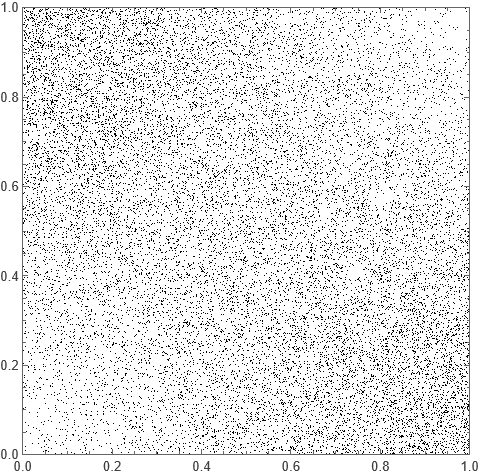}\\  \includegraphics[width=5cm]{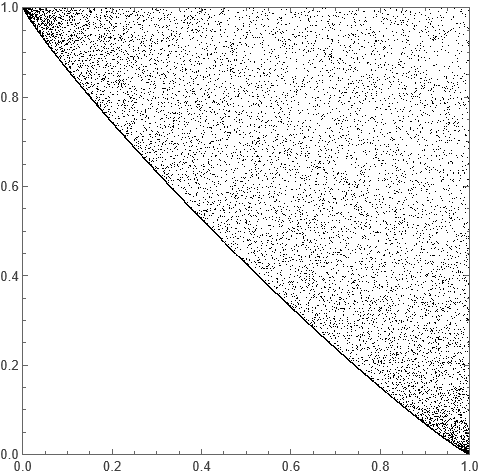}\ \ \includegraphics[width=5cm]{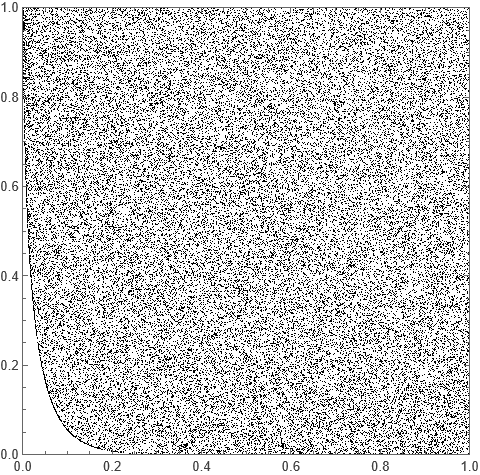}
\label{fig:1}
\end{figure}
\end{center}
  
  \begin{example}[Case of exponential shocks]\label{Ex:9}
    If 
    \[
        \begin{split}
           X\sim \exp(\lambda_1),\quad & Y\sim \exp(\lambda_2), \\
           Z_1\sim \exp(\mu_1),\quad  & Z_2\sim \exp(\mu_1),
        \end{split}
    \]
    then we have a copula for any  $\alpha,\beta\in(0,1]$ by letting
    \[
        C_{\alpha,\beta}(u,v)=\max\{0,uv-(u^\alpha-u)(v^\beta-v)\},\ \mbox{where}\ \alpha=\dfrac{\lambda_1} {\lambda_1+\mu_1},\beta=\dfrac{\lambda_2} {\lambda_2+\mu_2}.
    \]
  \end{example}
  
  \begin{proof}
    This can easily be determined by verifying conditions $(G_1)-(G_3)$ for the generating functions of the form
    \[
        f(t)=t^\alpha-t,\quad\mbox{for}\quad \alpha\in(0,1].
    \]
   More generally, a function
   \[
    f(t)=t^\alpha(1-t^\beta),
   \]
   is a generating function $f$ of an RMM copula if $\beta\in(0,1]$ for $\alpha=1$ and $\beta\geqslant 1-\alpha$ for $\alpha\in(0,1)$.
  \end{proof}

%%%%%%vstravi

%%%%%%od

\section{Survival copulas of RMM copulas}\label{sec:smm} 

For every RMM copula
\[
    C(u,v)=\max\{0,uv-f(u)g(v)\}
\]
we have its survival copula
\[
    \widehat{C}(u,v)=u+v-1+C(1-u,1-v)=\max\{u+v-1,uv-f(1-u)g(1-v)\}.
\]
We will now acquire the properties of functions $h,k:\II\rightarrow\II$ defined by $h(u)=f(1-u)$ and $k(v)=g(1-v)$; this time we need auxiliary functions
\[
 h_\dag(u)=\dfrac{h(u)}{1-u}, k_\dag(v)=\dfrac{k(v)}{1-v}, \widehat{h_\dag}(u)=u-h(u),\ \mbox{and}\ \  \widehat{k_\dag}(v)=v-k(v), 
\]
defined so that they have positive values. Here are their properties.
\begin{enumerate}[(S1)]
  \item $h(0)=k(0)=0,\ h(1)=k(1)=0$, and $\widehat{h_\dag}(1)=0,\widehat{k_\dag}(1)=0$;
  \item ${\widehat{h_\dag}}$ and ${\widehat{k_\dag}}$ are nondecreasing on $\II$; 
  \item functions $h_\dag(u)$ and $k_\dag(u)$ are nondecreasing on $\II$.
\end{enumerate}

\begin{lemma}
  If $f$ respectively $g$ has properties (G1)-(G3), then  $h$ respectively $k$ has properties (S1)-(S3).
\end{lemma}

\begin{proof}
  We consider only the case of $f$ and $h$ since the case of $g$ and $k$ goes similarly.
  Property (G1) clearly implies (S1). By (G2) function $\hat{f}(u)=f(u)+u$ is nondecreasing. So, $h(u)=f(1-u) =\widehat{f}(1-u)-1+u$ and $\widehat{h_\dag}(u)=u-h(u)=1-\widehat{f}(1-u)$
   is nondecreasing. Finally, since $f^*(u)=\dfrac{f(u)}{u}$ is nonincreasing, it follows that $h_\dag(u)=\dfrac{h(u)}{1-u}=\dfrac{f(1-u)}{1-u}=f^*(1-u)$ is nondecreasing.
\end{proof}

\begin{lemma}\label{lemma3}
  If functions $h$ and $k$ satisfy properties (S1)-(S3) then functions $f(u)=h(1-u)$ and $g(v)=k(1-v)$ satisfy properties (G1)-(G3).
\end{lemma}

\begin{proof}
  Property (S1) clearly implies (G1). By (S2) we have that $\widehat{h_\dag}(u)=u-h(u)$ is nondecreasing. Then $\widehat{f}(u)= u+f(u)= u+h(1-u)= 1-\widehat{h_\dag}(1-u)$ is nondecreasing as well. In the same way we prove that since $\widehat{k_\dag}$ is nondecreasing then $\widehat{g}$ is of the kind. Finally, (S3) implies that 
  \[
    f^*(u)=\frac{f(u)}{u}=\frac{h(1-u)}{u}=h_\dag(1-u),
  \]
  respectively $ g^*(v)=k_\dag(1-v)$ is nonincreasing.
\end{proof}

We call a copula of the form
\[
    C_{h,k}(u,v)=\max\{u+v-1, uv-h(u)k(v)\},
\]
where $h$ and $k$ satisfy conditions (S1)-(S3), a \emph{survival RMM copula} or \emph{SMM copula} for short. 

Let us now exhibit how the RMM and SMM copulas and their generating functions can be deduced from the original maxmin copula $C$ of Omladi\v{c} and Ru\v{z}i\'{c} \cite{OmRu}, and its generating functions. The RMM respectively SMM copula can be obtained as
\[
    C^{\sigma_2}(u,v)= u-C(u,1-v)\quad\mbox{respectively}\quad C^{\sigma_1}(u,v)= v-C(1-u,v).
\]
Choose generating functions $\phi$ and $\psi$ as in \cite{OmRu}, so that
\[
    C(u,v)=\min\{u,uv+(\phi(u)-u)(v-\psi(v)\}.
\]
Then
\[
     C^{\sigma_2}(u,v)=\max\{0,uv-(\phi(u)-u)(1-v-\psi(1-v)) \}
\]
and the generating functions of the so obtained RMM copula are clearly
\[
    f(u)=\phi(u)-u\quad\mbox{and}\quad g(v)=1-v-\psi(1-v).
\]
On the other hand,
\[
    C^{\sigma_1}(u,v)=\max\{u+v-1, uv-(\phi(1-u)-(1-u))(v-\psi(v))\},
\]
so that the generating functions of the so acquired SMM copula are

\[
    h(u)=\phi(1-u)-(1-u)\quad\mbox{and}\quad k(v)=v-\psi(v)).
\]

%%%%%%do tu

%%%%%%

\section{Stochastic model for SMM copulas}\label{sec:model}

In this Section we present a different model as an outcome of the stochastic setup exhibited in Section \ref{sec:setup}. As before, suppose that $(X,Y)$ and $(Z_1,Z_2)$ are two independent pairs of random variables, where $X$ and $Y$ are independent, while $Z_1$ and $Z_2$ are countermonotonic. More precisely, let $Z$ be a random variable distributed uniformly on $\II$ and let $Z_1=G_1^{-1}(Z)$ and  $Z_2=G_2^{-1}(I-Z)$ for some distribution functions $G_1$ and $G_2$. Unlike in Section \ref{sec:setup} we assume now that
\[
    U=\min\{X,Z_1\},\quad V=\min\{Y,Z_2\}.
\]

\begin{theorem}\label{thm7}
  The copula of pair $(U,V)$ is an SMM copula.
\end{theorem}

\begin{proof}
  Suppose that
\[
    C_{h,k}(u,v)=\max\{u+v-1, uv-h(u)k(v)\},
\]
where $h,k$ satisfy properties (S1)-(S3). By Lemma \ref{lemma3} functions $f(u)=h(1-u)$ and $g(u)=k(1-u)$ satisfy conditions (G1)-(G3). Define
\[
    C_{f,g}(u,v)=\max\{0,uv-f(u)g(v)\},
\]
then $C_{f,g}$ is an RMM copula by \cite[Theorem 3]{KoOm}. Let us compute its survival copula
\[
    \begin{split}
       \widehat{C_{f,g}}(u,v) & =u+v-1-C_{f,g}(1-u,1-v) \\
         & =\max\{ u+v-1, uv-f(1-u)g(1-v) \}.
    \end{split}
\]
According to \cite[Section 2.6]{Nels} this is a copula since it is a survival copula of a copula. We introduce $h(u)=f(1-u)$ and $k(v)=g(1-v)$ to conclude
\[
    \widehat{C_{f,g}}(u,v)=\max\{ u+v-1,uv-h(u)k(v)\}=C_{h.k}(u,v).
\]
\end{proof}

We combine Theorem \ref{thm8} with Theorem \ref{thm7} to obtain a shock model for an SMM copula. 

\begin{theorem}
  Let $U$ and $V$ be two random variables with respective distribution functions $F_U$ and $F_V$ such that there exists an $x_0\in\mathds{R}$ with $F_U(x_0),F_V(x_0)\in(0,1)$. Furthermore, let the copula of the pair $(F_U,F_V)$ be an SMM copula with generating functions $h$ and $k$, i.e.,
  \[
    C(u,v)=\max\{u+v-1,uv-h(u)k(v)\}.
  \]
  Then there exist random variables $X,Y,Z_1,Z_2$ such that $X,Y,Z_i$ are independent for $i=1,2$, $Z_2$ is countermonotonic to $Z_1$, and function
  \[
    H(x,y)=C(F_U(x),F_G(y))
  \]
  is the joint distribution function of the variables 
  \[\min\{X,Z_1\}\quad\mbox{and}\quad \min\{Y,Z_2\}.\]
\end{theorem}

\end{document}